\documentclass[11pt]{article}

\usepackage{amsmath,amsthm,amsfonts,amssymb,graphicx}

\def\Z{{\mathbb Z}}

\def\d{\Delta 334}

\def\mod{\text{ mod }}
\def\SLZ{SL_3(\Z)}
\def\SLTwo{SL_3(\Z/2\Z)}
\def\SLThree{SL_3(\Z/3\Z)}
\def\SLp{SL_3(\Z/p\Z)}

\newtheorem{theorem}{Theorem}[section]
\newtheorem{lemma}[theorem]{Lemma}
\newtheorem{conj}[theorem]{Conjecture}
\newtheorem{defn}[theorem]{Definition}

\title{The 334-Triangle Graph of $\SLZ$}
\author{Eric S. Egge\footnote{Corresponding author.} \\ Department of Mathematics and Statistics\\ Carleton College\\ 1 North College Street\\ Northfield, MN 55057 USA \\ eegge@carleton.edu\\ \\ Michaela A. Polley\\ Department of Mathematics and Statistics\\ Carleton College\\ 1 North College Street\\ Northfield, MN 55057 USA\\ polleym@carleton.edu}

\date{}

\begin{document}

\maketitle
\begin{abstract}
Long, Reid, and Thistlewaite have shown that some groups generated by representations of the $\d$ triangle group in $\SLZ$ are thin, while the status of others is unknown. 
In this paper we take a new approach: for each group we introduce a new graph that captures information about representations of $\d$ in the group.
We provide examples of our graph for a variety of groups, and we use information about the graph for $\SLTwo$ to show that the chromatic number of the graph for $\SLZ$ is at most eight.
By generating a portion of the graph for $\SLZ$ we show its chromatic number is at least four;  we conjecture it is equal to four.

\medskip

\end{abstract}

\section{Introduction}

Consider a subgroup $G$ of $SL_n(\Z)$, the group of all $n\times n$ matrices with integer entries and a determinant of one, under matrix multiplication. Later we will also consider subgroups of $SL_n(\Z/p\Z)$, the groups of all $n \times n$ matrices with determinant one and entries in $\Z/p\Z$, under matrix multiplication.

We say $G$ is a thin group whenever $G$ has infinite index in $SL_n(\Z)$ and the Zariski closure of $G$ is all of $SL_n.$ To define the Zariski closure of $G$, suppose 
$$A = \begin{pmatrix}
x_{11} && x_{12} && \dots && x_{1n}\\
x_{21} && x_{22} && \dots && x_{2n}\\
\vdots && \vdots && \ddots && \vdots\\
x_{n1} && x_{n2} && \dots && x_{nn}
\end{pmatrix}$$
is in $G$ and let $p$ be a polynomial in the $n^2$ variables $x_{11},x_{12},\dots,x_{1n},\dots,x_{n1},x_{n2},$ and $x_{nn}$. We define $p(A)$ to be $p$ evaluated at the entries of $A$. Let $I(G)$ be the set of all polynomials $p$ such that $p(A)=0$ for all $A\in G$. The Zariski closure of $G$ is the set of all matrices $B\in SL_n(\Z)$ such that for each polynomial $p\in I(G)$ we have $p(B)=0$. By construction, $G$ is a subset of its Zariski closure.

While groups that have both infinite index and are Zariski dense have been studied for the past 150 years, and they have been called thin for the last 15, there are still many groups for which we cannot definitively determine whether or not they are thin \cite{Kontorovich}. We explore this in three (non-)examples.

Consider the subgroup $G$ of $SL_2(\Z)$ generated by 
\begin{center}
    $A=$ $\left(\begin{matrix}
    1 && 1\\ 0 && 1
    \end{matrix}\right)$ and $B=$ $\left(\begin{matrix}
    0 && 1\\ -1 && 0
    \end{matrix}\right)$.
\end{center}
As Kontorovich, Long, Lubotzky, and Reid note \cite{Kontorovich}, it is well known that this group is all of $SL_2(\Z),$ and thus has index one and is not thin.

Next consider the subgroup $G$ of $\SLZ$ generated by
\begin{center}
    $A=$ $\left(\begin{matrix}
    0 && 0 && 1\\ 1 && 0 && 0 \\ 0 && 1 && 0
    \end{matrix}\right)$ and $B=$ $\left(\begin{matrix}
    1 && 2 && 3\\ 0 && -2 && -1 \\ 0 && 3 && 1
    \end{matrix}\right)$.
\end{center}
We can check that $G$ is a representation of the $\d$ triangle group, which is defined by
\begin{equation}
    T=\langle a,b\ |\ a^3=b^3=(ab)^4=e \rangle.
    \label{eq_triangle}
\end{equation}
In fact, this representation is faithful; one can use work of Margulis \cite{Margulis} to show that this implies $G$ has infinite index in $\SLZ$. It turns out that the Zariski closure of $G$ is $\SLZ$, so $G$ is thin.

By contrast, it is not always known whether a given subgroup of $SL_n(\Z)$ is thin. For example, consider the group $G$ generated by
\begin{center}
    $A=$ $\left(\begin{matrix}
    1 && 1 && 2\\ 0 && 1 && 1 \\ 0 && -3 && -2
    \end{matrix}\right)$ and $B=$ $\left(\begin{matrix}
    -2 && 0 && -1\\ -5 && 1 && -1 \\ 3 && 0 && 1
    \end{matrix}\right)$.
\end{center}
We can mod out by a prime number to show that this group is Zariski dense, as described in \cite{Kontorovich}, but it is not known whether the group has infinite index in $\SLZ$.

Many of the thin group candidates studied in \cite{SL3Z}, as well as the examples above, are representations of the $\d$ triangle group $T$ defined in \eqref{eq_triangle}. In fact, we have used ideas from \cite{SL3Z} to generate thousands of additional thin group candidates, all of which are also representations of $T$. In this paper we introduce and study a natural graph on the set of elements of order three in a group $G$ that captures interesting information about the set of representations of $\d$ in $G$. 

For any group $G$, let $\d(G)$ be the graph whose vertices are the elements $a\in G$ such that $a^3=e$, in which there is an edge between two vertices $a$ and $b$ if and only if $(ab)^4=e$. 
We note that since $(ab)^4 = e$ defines a symmetric relation, $\d(G)$ is an undirected graph.
We call $\d(G)$ the {\em 334-triangle graph of $G$}.
In this paper we explore the properties of $\d(G)$ generally before considering a number of finite examples. 
We then turn our attention to $\d(\SLZ)$.

Although we do not know of an explicit connection between $\d(\SLZ)$ and thin groups in $\SLZ$, this graph seems to be of independent interest.
Among the properties of $\d(\SLZ)$ we could study, we focus on the chromatic number.
We use the natural homomorphism from $\SLZ$ to $\SLTwo$ to show that this chromatic number is at most eight, and we examine a small portion of the graph to show it is at least four.
We conjecture it is equal to four.

\section{The 334-Triangle Graph}

In this section, we analyze the 334-triangle graph for a variety of groups. To start, we prove four facts about $\d(G)$: the identity element is adjacent to itself and nothing else, the identity element is the only element that is adjacent to itself, every element is adjacent to its inverse, and for Abelian groups these are the only edges. We also describe $\d(G)$ when $G$ is a direct product of two groups.

\begin{lemma}
For any group $G$, the identity in $G$ is adjacent in $\d(G)$ to itself and is adjacent to no other vertex in $\d(G)$. Furthermore, the identity is the only element in $\d(G)$ that is adjacent to itself.
\label{Lm_Identity}
\end{lemma}

\begin{proof}
Let $G$ be a group and let $e$ be the identity element of $G$. Then $e^3=e$ and $(ee)^4=e$, so $e$ is a vertex in $\d(G)$ and is adjacent to itself. 

Let $A$ be an order three element of $G$, so that $A^3=e,$ but $A\ne e$. Then, $(Ae)^4=A$, so there is no edge connecting $A$ to $e$. Hence, $e$ is not adjacent to any other element in the $\d$ graph of $G$. 

Finally, suppose $A$ is adjacent to itself. Since $A$ is a vertex in $\d(G)$, we have $A^3=e$ and therefore $A^9=e$. Since $A$ is adjacent to itself, we also have $(A^2)^4=e$, so $A^8=e$. Since $A^9=A^8,$ we must have $A=e$.
\end{proof}

In view of Lemma \ref{Lm_Identity}, we will almost always disregard the identity vertex in further discussions of $\d(G)$ and focus only on the non-identity component(s).

\begin{lemma}
For any group $G$, and any vertex $A$ in $\d(G)$, the element $A^{-1}$ is also a vertex in $\d(G)$, and $A$ and $A^{-1}$ are adjacent.
\label{Lm_Inverses}
\end{lemma}

\begin{proof}
Let $G$ be a group and let $A\in G$ be an element of $G$ such that $A^3=e$. If $A=e$, then $A=A^{-1}$, and by Lemma \ref{Lm_Identity}, $A$ is adjacent to itself. Therefore, $A$ and $A^{-1}$ are both vertices in $\d(G)$ and are adjacent. 

Now suppose $A\ne e$. Since $|A|=|A^{-1}|=1$ and $A^3=e$, we also know that $(A^{-1})^3=e$. Therefore, $A^{-1}$ is a vertex in $\d(G)$. In addition, $(AA^{-1})^4=e^4=e$, so $A$ and $A^{-1}$ are adjacent.
\end{proof}

Using Lemmas \ref{Lm_Identity} and \ref{Lm_Inverses}, we are able to describe $\d(G)$ completely when $G$ is Abelian.  

\begin{lemma}
\label{lem:Abelian}
For any Abelian group $G$, two vertices $A,B$ in $\d(G)$ are adjacent if and only if $A^{-1}=B$.
\label{Lm_Abelian}
\end{lemma}

\begin{proof}
Let $G$ be an Abelian group and let $A$ and $B$ be two elements of $G$ such that $A^3=B^3=e$. We know $A$ and $B$ will be adjacent if and only if $(AB)^4=e$. However, $(AB)^4=A^4B^4=AB$. Thus, $A$ and $B$ will be adjacent if and only if $AB=e$, and this is only true when $A^{-1}=B$.
\end{proof}

We know that in some cases the converse of Lemma \ref{lem:Abelian} can fail. That is, there are some non-Abelian groups that also have the property that two vertices $A,B$ are adjacent if and only if $A^{-1}=B$. For example, this will happen if $G$ is not Abelian and has odd order divisible by three. It is not known whether the converse holds for groups whose order is divisible by six and that are generated by their order two and order three elements.

In the case where $H\subseteq G$, it turns out that $\d(H)$ is a subgraph of $\d(G)$. That is, $\d(H)$ is a subset of the vertices in $\d(G)$ along with all edges connecting them. Another interesting question is for $H\triangleleft G$: how are $\d(G)$, $\d(H)$ and $\d(G/H)$ related? This is a question for further research.

Finally, we also consider the 334-triangle graph of a direct sum of groups, showing that the graph $\d(G \oplus H)$ is the Kronecker product of $\d(G)$ and $\d(H)$.

\begin{defn}
For any graphs $G$ and $H$ with vertices $g,g'\in G$ and $h,h'\in H$, the Kronecker product of $G$ and $H$, written $G \otimes H$, is the graph with vertices of the form $(g,h)$ and in which $(g,h)$ and $(g',h')$ are adjacent if and only if both $g$ and $g'$ are adjacent in $G$ and $h$ and $h'$ are adjacent in $H$.
\label{def_kronecker}
\end{defn}

\begin{lemma}
For any groups $G$ and $H$, we have $\d(G\oplus H)=\d(G)\otimes\d(H).$
\label{Lm_KronProduct}
\end{lemma}

\begin{proof}
Let $G$ be a group with elements $g$ and $g'$ such that $g^3=g'^3=e_G$ and let $H$ be a group with elements $h$ and $h'$ such that $h^3=h'^3=e_H,$ where $e_G$ and $e_H$ are the identity elements of $G$ and $H$, respectively. Consider elements $(g,h)$ and $(g',h')$ in $G\oplus H$. We wish to show that $(g,h)$ and $(g',h')$ will be adjacent in $\d(G\oplus H)$ if and only if $g$ and $g'$ are adjacent in $\d(G)$ and $h$ and $h'$ are adjacent in the $\d(H)$.

In order for an element $(g,h)$ in $G\oplus H$ to have order three, both $g$ and $h$ must have either order three or one. Thus, $g$ and $h$ will both appear in $\d(G)$ and $\d(H)$, respectively. If $g$ and $g'$ are adjacent in $\d(G)$ and $h$ and $h'$ are adjacent in $\d(H)$, then $((g,h)(g',h'))^4=((gg')^4,(hh'^4))=(e_G,e_H)$, so $(g,h)$ and $(g',h')$ are adjacent in $\d(G \oplus H)$

On the other hand, let $(g,h)$ and $(g',h')$ be adjacent in $\d(G\otimes H)$. Then, $((g,h)(g',h'))^4=((gg')^4,(hh'^4))=(e_G,e_H)$, so $(gg')^4=e_G$ and $(hh')^4=e_H$, so $g$ and $g'$ are adjacent in $\d(G)$ and $h$ and $h'$ are adjacent in the $\d(H)$.
\end{proof}

Having now laid the groundwork for understanding these graphs, and having completely described the graph for all Abelian groups and groups that are isomorphic to direct sums of two or more groups, let us consider $\d(G)$ of some groups that are neither of these. 
We will begin by analyzing $S_4$, the set of permutations on four letters. $S_4$ has eight order three elements: $(123),(132),(124),(142),(134),(143),(234),$ and $(243)$. 
Thus, the non-identity component of the graph will have eight vertices. 
We are interested in the subgroup these elements generate, since any element of the group that is not part of our subgroup will not appear in the graph. 
For $S_4$, the relevant subgroup is $A_4$. 
Therefore, the graphs for $S_4$ and $A_4$ are the same, so we will focus our discussion on $A_4$. 
(Indeed, the same reasoning shows $\d(S_n) = \d(A_n)$ for all $n \ge 2$.)

We can check that the order three elements of $A_4$ fall into two conjugacy classes, each with four elements.
As we will show, each order three element is adjacent to exactly those order three elements to which it is not conjugate.
To see this, note that for two order three elements in $S_4$ there are two possibilities for their product:
\begin{equation}
    (abc)(abd)=(ac)(bd)
    \label{Eq_S4yes}
\end{equation}
and
\begin{equation}
    (abc)(bad)=(adc).
    \label{Eq_S4No}
\end{equation}
Since these are the only two options, and $(adc)^4\ne e$, there will be an edge connecting elements $A$ and $B$ if and only  if $A$ is of the form $(abc)$ and $B$ is of the form $(abd)$. 
We can check that this occurs exactly when $A$ and $B$ are not conjugate in $A_4$.

\begin{figure}[htbp]
    \centering
    \includegraphics{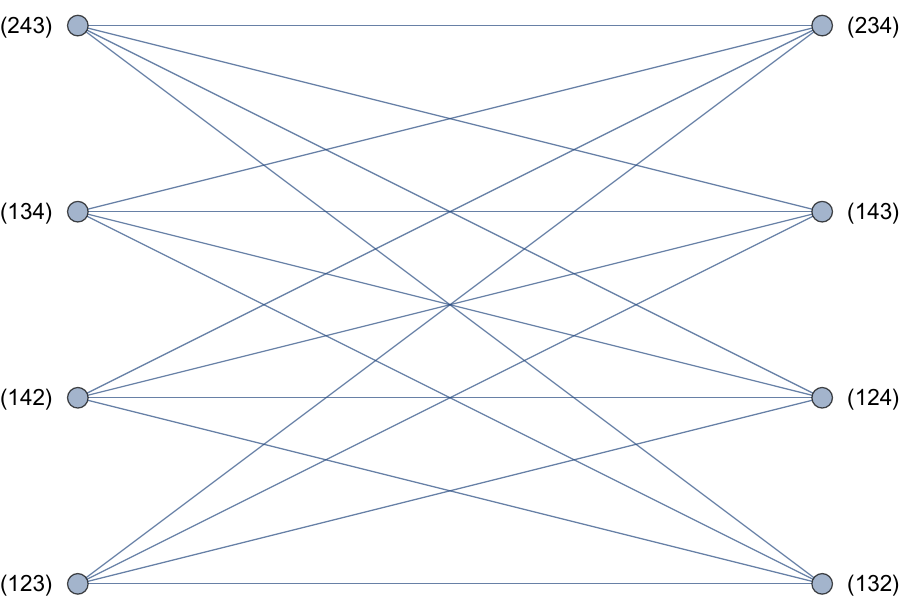}
    \caption{The non-identity component of $\Delta334(A_4)$.}
    \label{fig_A_4}
\end{figure}

Based on the analysis above, we find that the graph $\d(A_4)$ is as in Figure \ref{fig_A_4}.
This is a complete bipartite graph, with the conjugacy classes forming the bipartition. This graph has chromatic number two and has cycles of length four, six, and eight; and its clique number is two. 

We leave it to the reader to verify the interesting fact that $SL_2(\Z/3\Z)$ and $S_4$, which have the same order, also have isomorphic 334-triangle graphs.

\begin{figure}[htbp]
    \centering
    \includegraphics{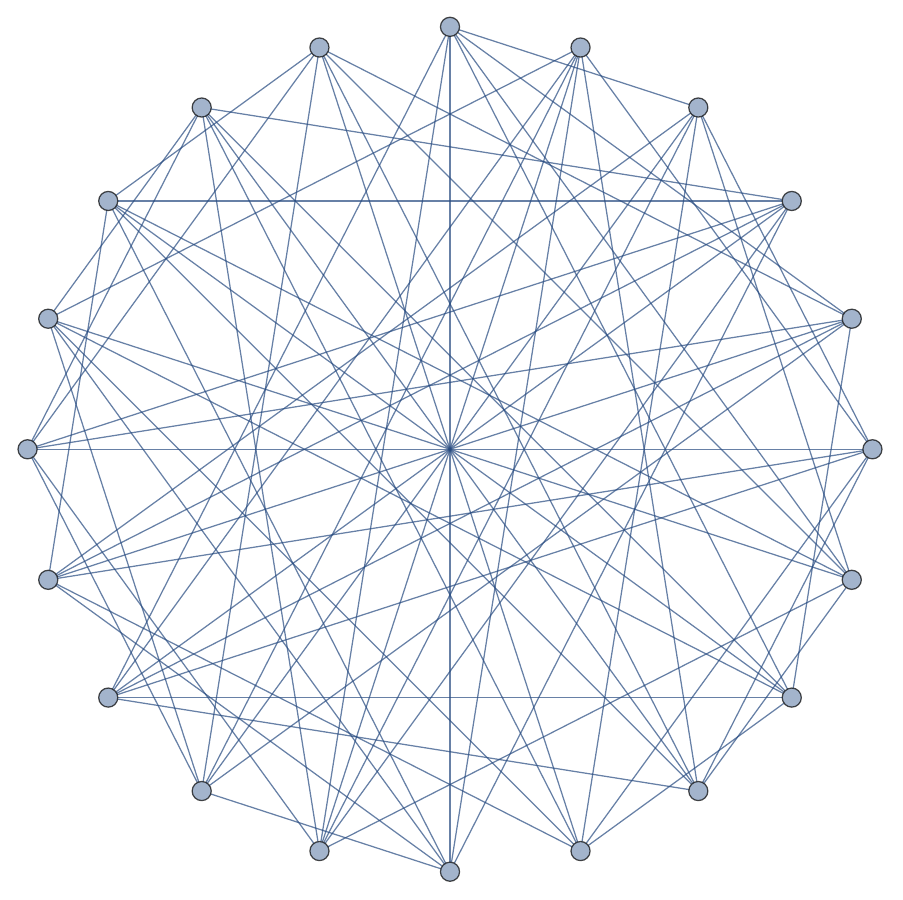}
    \caption{The non-identity component of $\Delta334(S_5)$.}
    \label{fig_S5}
\end{figure}

As another example, we observe that the 334-triangle graph of $S_5$ (and therefore of $A_5$) is not bipartitie -- $(125),(124)$, and $(123)$ form a three-cycle -- and every vertex has degree seven. This graph is shown in Figure \ref{fig_S5}.

We are most interested in $\Delta 334(\SLZ)$, since we have generated many candidates for thin groups that are subgroups of $\SLZ$, as mentioned in the Introduction. However, this graph is infinite, so we will begin with finite graphs that are images of this graph after modding out by a prime. First we analyze the non-identity component of $\Delta 334(\SLTwo)$, which is shown in Figure \ref{fig_SL(3,2)}. We used \textit{Mathematica} to generate all of the order three elements of $SL_3(\Z/2\Z)$. There are $56$ such elements, thus, the non-identity component of $\d(\SLTwo)$ contains $56$ vertices. Every vertex has degree $19$ and there are cycles of all lengths from $3$ to $56$. In particular, this component is connected and Hamiltonian. It also has chromatic number eight and cliques up to size five. One such clique is 

$$\left\{\begin{pmatrix}0 && 1 && 0\\1 && 0 && 1\\1 && 1 && 0\end{pmatrix},
\begin{pmatrix}1 && 0 && 1\\0 && 1 && 1\\0 && 1 && 0\end{pmatrix},
\begin{pmatrix}1 && 1 && 1\\1 && 0 && 0\\0 && 0 && 1\end{pmatrix},
\begin{pmatrix}1 && 0 && 0\\1 && 1 && 1\\1 && 1 && 0\end{pmatrix},
\begin{pmatrix}1 && 0 && 0\\1 && 0 && 1\\0 && 1 && 1\end{pmatrix}
\right\}.$$

\begin{figure}[htbp]
    \centering
    \includegraphics{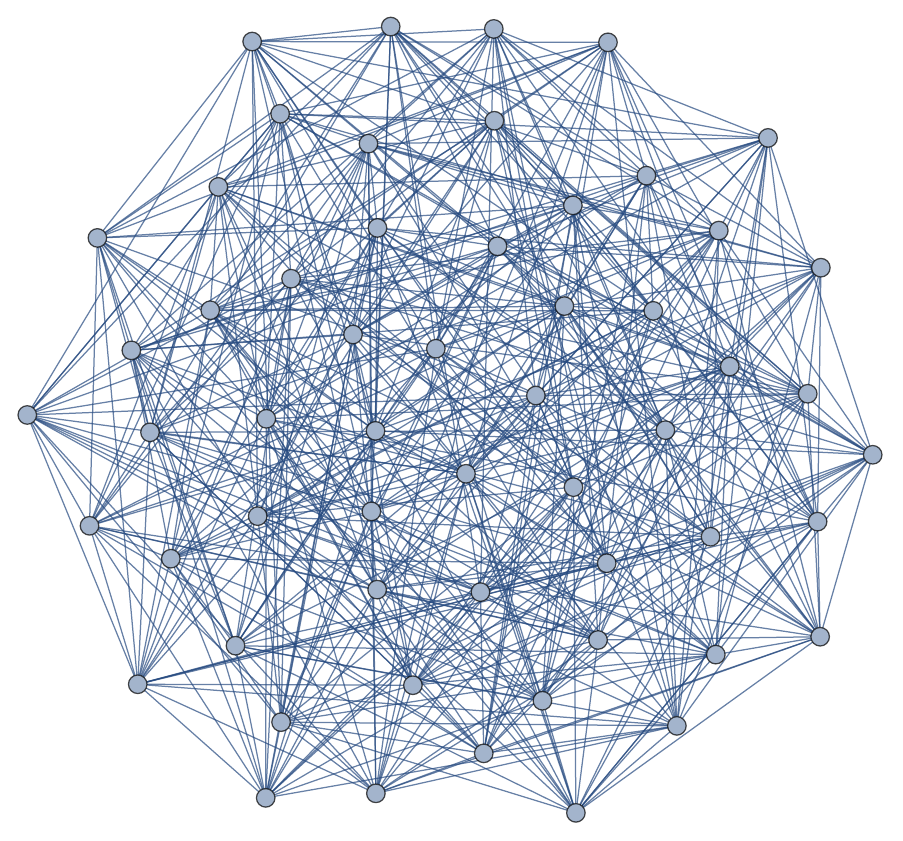}
    \caption{The non-identity component of $\d(\SLTwo)$.}
    \label{fig_SL(3,2)}
\end{figure}
We are able to say less about $\d(\SLThree)$ because the graph becomes so much larger. We know there are $728$ matrices of order three in $\SLThree$, so the non-identity component of $\d(\SLThree)$ will contain $728$ vertices. We know there is a single connected non-identity component. All vertices have degree either 118 or 136. We were unable to determine the length of cycles or chromatic number of this graph due to its size.

\section{$\d(\SLZ)$}

We now consider $\d(\SLZ)$. 
We know $\SLZ$ is infinite, and we can show it has an infinite number of order three elements. For example, for any $a,b,c\in \Z$, 
$$\begin{pmatrix}
1&& 3a && 3b\\
0 && -2-3c && -1-3c-3c^2\\
0 && 3 && 1+3c
\end{pmatrix}$$
has integer entries, determinant one, and order three. Thus, $\d(\SLZ)$ is also infinite. However, using what we know about $\d(\SLTwo)$ we are able to put bounds on the chromatic number of $\d(\SLZ)$. 
We will do this by showing that any edge in $\d(\SLZ)$ reduces to an edge in $\d(\SLTwo)$. 
We will then lift a proper coloring from $\d(\SLTwo)$ to $\d(\SLZ)$ by coloring each vertex in $\d(\SLZ)$ the same color as its image in $\d(\SLTwo)$. 

Before we begin, we note that we are using the fact that the natural homomorphism from $\SLZ$ to $\SLTwo$ induces a graph homomorphism from $\d(\SLZ)$ to $\d(\SLTwo)$.
It's worth noting that in general any homomorphism from a group $G$ to a group $H$ will similarly induce a graph homomorphism from $\d(G)$ to $\d(H)$, underscoring the fact that the $334$-Triangle graph is natural from an algebraic point of view.

We first show that no vertex in $\d(\SLZ)$ reduces to the identity modulo any prime.

\begin{lemma}
For any matrix $A\in \SLZ$ with order three, $A$ cannot be equivalent to the identity matrix modulo any prime.
\label{Lm_ReduceToI3}
\end{lemma}

\begin{proof}
Let $A$ be a matrix in $SL_3(\Z)$ with $A^3=I_3$, where $I_3$ is the identity element of $SL_3(\Z)$, and assume $A \neq I_3$.
Let $p$ be prime, and assume by way of contradiction that $A$ is equivalent to the identity modulo $p$.

Since $A^3 = I_3$, the minimal polynomial for $A$ over ${\mathbb Q}$ must divide $x^3-1 = (x-1)(x^2+x+1)$.
Since $A \neq I_3$, the minimal polynomial cannot be $x-1$.
We also note that $x^2+x+1 = (x-\omega)(x-\overline{\omega})$, where $\omega = -\frac12 + \frac{\sqrt{3}}{2} i$ is a primitive third root of unity.
Since $\omega$ and $\overline{\omega}$ are not rational, the minimal polynomial must be $x^2+x+1$ or $x^3-1$.
But if the minimal polynomial were $x^2+x+1$ then by the Cayley-Hamilton theorem the characteristic polynomial would be $-(x-\omega)(x^2+x+1)$ or $-(x-\overline{\omega})(x^2+x+1)$, neither of which has constant term equal to $\det A = 1$.
Therefore, the minimal polynomial for $A$ over $\mathbb Q$ is $x^3-1$.

The fact that the minimal polynomial for $A$ over $\mathbb Q$ is $x^3-1$ implies $A$ has an eigenvector $\vec{v}_1$ with eigenvalue one, and that there is a vector $\vec{v}_2 \neq \vec{0}$ with $(A^2+A+I_3)\vec{v}_2 = \vec{0}$.
We set $\vec{v}_3 = A \vec{v}_2$.
By scaling if necessary, we can assume all of the entries of $\vec{v}_1$, $\vec{v}_2$, and $\vec{v}_3$ are integers and that the entries of $\vec{v}_2$ have no common prime factor.
We can check that $A\vec{v}_3 = -\vec{v}_2 - \vec{v}_3$ and that $\vec{v}_1$, $\vec{v}_2$, and $\vec{v}_3$ form a basis for ${\mathbb Q}^3$.
Therefore, there is an invertible $3 \times 3$ matrix $M$ with entries in $\mathbb Z$ such that
\[ A M = M \left(\begin{matrix} 1 & 0 & 0 \\ 0 & 0 & -1 \\ 0 & 1 & -1\end{matrix}\right).\]
In particular, the columns of $M$ are $\vec{v}_1$, $\vec{v}_2$, and $\vec{v}_3$.
We can now use our assumption that $A$ is equivalent to the identity modulo $p$ to check that the entries of the middle column of $M$ are all divisible by $p$.
But this contradicts the fact that the entries of $\vec{v}_2$ do not have a common prime factor.
\end{proof}

Having shown that no matrices in the non-identity component of $\d(\SLZ)$ reduce to the identity modulo a prime $p$, we can show that every edge in $\d(\SLZ)$ maps to an edge in $\d(\SLp)$.

\begin{lemma}
Let $A$ and $B$ with $A\ne B$ be adjacent vertices in $\d(\SLZ)$. Then, for $A'=A \mod p$ and $B'=B \mod p$, $A'$ and $B'$ are adjacent in $\d(\SLp)$ and $A'\ne B'$.
\label{Lm_mappingEdges}
\end{lemma}
\begin{proof}
Let $A,B$ be adjacent vertices in $\d(\SLZ)$ with $A\ne B$ and let $p$ be prime. By Lemma \ref{Lm_Identity}, neither $A$ nor $B$ are the identity, and by Lemma \ref{Lm_ReduceToI3} neither reduce to the identity modulo $p$. 

Let $C=AB$ and let $C'=A'B'$. We know that $C'=C \mod p$ by the rules of modular arithmetic. We know that $C^4=I_3$, so $C'^4$ must also equal the identity. Thus, $A'$ and $B'$ are adjacent in $\d(\SLp)$. By Lemma \ref{Lm_Identity}, the only loop in $\d(\SLp)$ is at the identity and $A\ne B\ne I_3$, so $A'\ne B'$.
\end{proof} 

Finally, since every edge in $\d(\SLZ)$ maps to an edge in $\d(\SLp)$, we are able to set an upper bound on the chromatic number, $\chi(\d(\SLZ)$.

\begin{theorem}
For any prime $p$, we have $\chi(\d(\SLZ))\leq \chi(\d(\SLp))$.
\label{Lm_upperbound}
\end{theorem}
\begin{proof}
If we have a proper coloring of $\d(\SLp)$ where $p$ is prime, then we can lift it to a proper coloring of $\d(\SLZ)$. We do this by coloring every vertex in $\d(\SLZ)$ the same color as its image in $\d(\SLp)$. For any two adjacent vertices in $\d(\SLZ)$, their images are also adjacent in $\d(\SLp),$ so they will have different colors. Thus, we will have a proper coloring of $\d(\SLZ).$ Therefore, the chromatic number of $\d(\SLZ)$ is at most the chromatic number of $\d(\SLp)$.
\end{proof}

Since $\chi(\d(\SLTwo))=8$, and two is prime, by Theorem \ref{Lm_upperbound} we know that the chromatic number of $\d(\SLZ)$ is at most eight. We have used \textit{Mathematica} to generate a finite portion of this graph with about 25,000 vertices; the chromatic number of this portion is four. Thus, one lower bound for the chromatic number of $\d(\SLZ)$ is four. We conjecture that $\chi(\d(\SLZ))$ is exactly four.
\begin{conj}
$\chi(\d(\SLZ))=4$.
\end{conj}
Using the IGraph/M package, which uses the Boyer–Myrvold algorithm, we have found that the  non-identity part of $\d(\SLZ)$ that we have generated so far is non-planar. Thus, the overall graph is also non-planar. Additionally, it is connected, has cycles of varying lengths, and has no cliques of size greater than three. We conjecture that all three of these facts hold for the entire graph $\d(\SLZ)$.
\begin{conj}
The non-identity component of $\d(\SLZ)$ is connected.
\end{conj}
\begin{conj}
There are cycles of every length in $\d(\SLZ)$.
\end{conj}
\begin{conj}
There are no cliques of size greater than three in $\d(\SLZ)$.
\end{conj}

\section{Acknowledgements}

The authors would like to thank Gabriel Konar-Steenberg for his support at the start of this project writing Python code that allowed us to explore this question.

\bibliographystyle{plain}
\bibliography{references}

\end{document}